\documentclass[showpacs,showkeys,amssymb,aps]{revtex4}

\usepackage{amsmath}

\usepackage{graphicx}

\begin{document}

\author{Richard J. Mathar}
\email{mathar@strw.leidenuniv.nl}
\homepage{http://www.strw.leidenuniv.nl/~mathar}
\affiliation{
Leiden Observatory, Leiden University, P.O. Box 9513, 2300 RA Leiden, The Netherlands}

\date{\today}
\title{Two Integrals of Geodetic Lines in Oblate Ellipsoidal Coordinates}

\begin{abstract}
The manuscript establishes a series expansion of the core integral that
relates changes in longitude and latitude along the geodetic line in
oblate elliptical coordinates, and of a companion integral which is the
path length along this line as a function of latitude.
The expansion is a power series in the scaled (constant) altitude
of the trajectory over the surface of the ellipsoid.
Each term of this series is reduced to sums over
inverse trigonometric functions, square roots and Elliptic Integrals.
The aim is to avoid purely
numerical means of integration.
\end{abstract}

\maketitle

\section{Scope}
\subsection{Geodetic Coordinates}

An ellipsoid is a reference surface fixed
by an equatorial radius $\rho_e$ and a polar radius $\rho_e$. In many applications
the second eccentricity $e$,
\begin{equation}
\rho_p^2= \rho_e^2(1-e^2),
\end{equation}
is the principal reduced parameter.
The three-dimensional ellipsoidal coordinates, altitude $h$ and the angles
of longitude $\lambda$ and latitude $\phi$, are basically
defined with the aid of a straight plumb line along the shortest distance
between a general point and its ``foot'' point on the surface.
The relation between
the Cartesian geocentric coordinates $(x,y,z)$ and the curvilinear $(h,\lambda,\phi)$ is
\cite{FukushimaJG79,VermeilleJG78,JonesJG76,PollardJG76,ZhangJG79,HradilekBullG50,KeelerSIR40},
\begin{equation}
\left(
\begin{array}{c}
x \\
y \\
z \\
\end{array}
\right)
=
\left(
\begin{array}{c}
\left[ N(\phi)+h \right] \cos \phi \cos \lambda \\
\left[ N(\phi)+h \right] \cos \phi \sin \lambda \\
\left[ N(\phi) (1-e^2)+h \right]\sin\phi \\
\end{array}
\right),
\end{equation}
where
\begin{equation}
N(\phi)\equiv \frac{\rho_e}{\sqrt{1-e^2\sin^2\phi}}
\label{eq.Ndef}
\end{equation}
is the distance between the foot point and the polar axis measured along
the straight extension of the plumb line.

The projection $\tau$ onto the polar axis,
\begin{equation}
\tau\equiv \sin\phi,
\end{equation}
will be useful to substitute trigonometric functions by rational functions.

\subsection{Geodetic Lines}

A geodetic line is the line of shortest Euclidean distance between two points
within a surface of constant height $h$.
This balances differential changes in the trajectory $\phi(\lambda)$ with
the two principal curvatures at each point; 
in consequence, the meridional radius of curvature
\begin{equation}
M(\tau) = \frac{\rho_e(1-e^2)}{(1-e^2\tau^2)^{3/2}}
\end{equation}
often appears to condense the notation.

The solution of the differential equations of the geodetic lines 
crystallizes in the integral \cite{MatharArxiv0711}
\begin{equation}
I(\tau)\equiv
\int d\tau
\frac
{c(h+M)}
{(N+h)^2(1-\tau^2)\sqrt{1-\tau^2-\frac{c^2}{(N+h)^2}}}
=\Delta \lambda
,
\label{eq.lint}
\end{equation}
which relates a difference in latitude---the limits of the integral---to
a difference in longitude---the right hand side. The obliquity parameter $c$ picks
an individual geodetic line out of the bundle of all lines that cross
a general point. Considering the $\tau$ at which the discriminant
of the root in the denominator of (\ref{eq.lint}) is zero shows that
$c$ is also the distance to the polar axis at the point highest above (or below)
the equatorial plane \cite{MatharArxiv0711}.

The single interest of this manuscript
is in demonstrating a semi-numerical approach to evaluation of this integral.
The constant of integration is tacitly fixed to imply the lower limit $\tau=\phi=0$
in the integral, because such a reference to the nodal line leads to well defined
branch cuts of all square roots involved.
The strategy is to expand the integrand into a power series of $h/\rho_e$,
and to exchange the order of integration and summation, which defines a family of integrals
with two additional parameters reminiscent of the order in the expansion.
Each of these is reduced to the level of multiple---but finite---sums over
Elliptic Integrals, assuming that these are accessible through a numerical library
\cite{CarlsonNumAlg10,CarlsonTOMS7,BulirschNumMath7,BulirschNumMath13}.

In overview, one way of computation of (\ref{eq.lint})
is addressed in Section \ref{sec.llc}.
Auxiliary integrals fall into two classes, one reducible to roots
and inverse trigonometric functions, the other to elliptic integrals.
The distance along the geodetic line defines another integral which
is treated in the same spirit in Section \ref{sec.sint}. Its power
series yields a family of integrals which can be recast for efficient
reuse of the functionality build in Section \ref{sec.llc}.
Finding the inverse of $I$ with respect to the parameter $c$ is closely
related to the inverse problem of geodesy and shortly addressed
in Section \ref{sec.inv}.

\section{Longitude-Latitude Coupling Integral}\label{sec.llc}
\subsection{Taylor Expansion in Powers of Altitude}
\begin{figure}
\includegraphics[width=12cm]{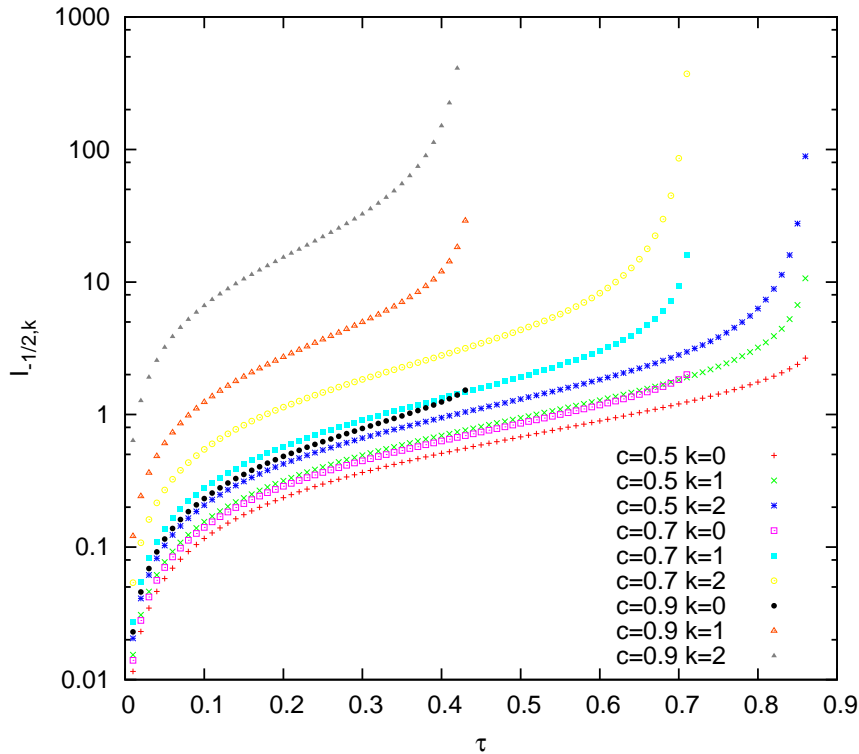}
\caption{Values of $I_{-1/2,k}$ with a lower limit of $\tau=0$ for three different
$c$ and three different $k$ with $e=0.08182$. Small $c$ indicates near-polar routes
where the branch cut with $I_{\beta,k}\to \infty$ is at larger $\tau$.}
\label{fig.1}
\end{figure}

Expansion of the auxiliary $N$ and $M$ and lifting of some square roots
provides a long write-up of (\ref{eq.lint}),
\begin{eqnarray}
I
&=&
\frac{c}{\rho_e}
\int d\tau
\frac
{\frac{h}{\rho_e}(1-e^2\tau^2)^{3/2}+1-e^2}
{(1+\frac{h}{\rho_e}\sqrt{1-e^2\tau^2})(1-\tau^2)
\sqrt{1-e^2\tau^2}
\sqrt{ (1+\frac{h}{\rho_e}\sqrt{1-e^2\tau^2})^2(1-\tau^2)-(c/\rho_e)^2(1-e^2\tau^2) }}
,
\label{eq.I}
\end{eqnarray}
which we intend to calculate.
The altitude $h$ and parameter $c$ appear only scaled with $\rho_e$, so
introducing a function of two dimensionless variables $h$ and $c$,
\begin{equation}
I_\alpha(h,c)
\equiv
\int d\tau
\frac
{(1-e^2\tau^2)^{\alpha}}
{(1+h\sqrt{1-e^2\tau^2})(1-\tau^2)
\sqrt{ (1+h\sqrt{1-e^2\tau^2})^2(1-\tau^2)-c^2(1-e^2\tau^2) }}
,
\label{eq.Ipows}
\end{equation}
shows the composition
\begin{equation}
I = \frac{c}{\rho_e}\left[
\frac{h}{\rho_e} I_1(\frac{h}{\rho_e},\frac{c}{\rho_e})+ (1-e^2)I_{-1/2}(\frac{h}{\rho_e},\frac{c}{\rho_e})
\right] 
.
\label{eq.IofIalph}
\end{equation}
The structure of the integrand is dominated by the two variables
\begin{equation}
T\equiv 1-\tau^2; \quad E\equiv 1-e^2\tau^2.
\label{eq.TEdef}
\end{equation}
The power series of (\ref{eq.Ipows}) becomes
\begin{equation}
I_\alpha(h,c)
= \sum_{s=0}^\infty (-h)^s \sum_{k=0}^s \kappa_{s,k} \int d\tau \frac{ E^{\alpha+s/2}}{(T-c^2E)^{k+1/2}} T^{k-1}
.
\label{eq.Ialphser}
\end{equation}
The coefficients $\kappa$ emerge from the product of the
geometric series of $1/(1+hE^{1/2})$ by the binomial
expansion of $1/\sqrt{(1+hE^{1/2})^2T-c^2E}$ in
(\ref{eq.Ipows}):
\begin{equation}
\kappa_{s,k} =
4^k \binom{-1/2}{k}
\sum_{l=k}^{\min(2k,s)}
\binom{k}{l-k}
(-1/2)^l
=
(-1)^k\binom{-1/2}{k} + (-1)^s 2^{2k-s-1} \binom{-1/2}{k} P_{2k-s-1}^{(1+s-k,-k)}(0)
.
\label{eq.kappa}
\end{equation}
The term with the Jacobi Polynomial $P$ is to be interpreted as zero if $k \le s/2$.
\begin{table}
\begin{tabular}{r|rrrrrrrrrrrr}
$s\backslash k$ & 0 & 1 & 2 & 3 & 4 & 5 & 6 & 7 & 8 & 9\\
\hline
0 & 1 & \\
1 & 1 & 1 & \\
2 & 1 & 1/2 & 3/2 & \\
3 & 1 & 1/2 & 0 & 5/2 & \\
4 & 1 & 1/2 & 3/8 & -5/4 & 35/8 & \\
5 & 1 & 1/2 & 3/8 & 5/8 & -35/8 & 63/8 & \\
6 & 1 & 1/2 & 3/8 & 5/16 & 35/16 & -189/16 & 231/16 & \\
7 & 1 & 1/2 & 3/8 & 5/16 & 0 & 63/8 & -231/8 & 429/16 & \\
8 & 1 & 1/2 & 3/8 & 5/16 & 35/128 & -63/32 & 1617/64 & -2145/32 & 6435/128 & \\
9 & 1 & 1/2 & 3/8 & 5/16 & 35/128 & 63/128 & -693/64 & 4719/64 & -19305/128 & 12155/128 & \\
\end{tabular}
\caption{
Table of the rational values of $\kappa_{s,k}$ by equation (\ref{eq.kappa}).
Each column attains a constant value for rows $s\ge 2k$.
}
\label{tab.kappa}
\end{table}
(\ref{eq.Ialphser}) states the problem in terms of integrals
\begin{equation}
I_{\beta,k}
\equiv \int d\tau \frac {T^{k-1} E^\beta} {(T-c^2E)^{k+1/2}}
=
\int d\tau \frac{(1-\tau^2)^{k-1} (1-e^2\tau^2)^\beta } { [1-\tau^2-c^2(1-e^2\tau^2)]^{k+1/2} }
.
\label{eq.Ibetakdef}
\end{equation}
for $\beta=\alpha+s/2=-1/2,0,1/2,1,3/2,\ldots$ and $k=0,1,2,3,\ldots$.
For small eccentricities, the $I_{\beta,k}$ are nearly independent of $\beta$
because $E$ is close to unity,
so it suffices to illustrate the values for zero lower limit
and variable upper limit $\tau$ for one value of $\beta$ in Figure \ref{fig.1}.
The main disadvantage of the method is that the eventual constancy of $\kappa_{s,k}$
down the columns of Table \ref{tab.kappa}
in conjunction with the alternating sign of $(-h)^s$ in (\ref{eq.Ialphser})
induces oscillatory behavior (cancellation effects) of the series if $h$ is not small.

For integrals classifications, (\ref{eq.Ibetakdef}) is phrased as
\begin{equation}
I_{\beta,k}
=
\frac{e^{2\beta}}{(1-c^2e^2)^{k+1/2}}
\int_0 d\tau \frac{(1-\tau^2)^{k-1} (a^2-\tau^2)^\beta }
  { (b^2-\tau^2)^{k+1/2} }
,
\label{eq.Iofab}
\end{equation}
where two parameters $a$ and $b$,
\begin{equation}
a\equiv 1/e^2 >  b\equiv (1-c^2)/(1-c^2e^2)
,
\label{eq.adef}
\end{equation}
decide on the branches of the Elliptic Integrals.
The exponent $\beta$ either leads to reduction to elementary functions
if it is an integer (Section \ref{sec.elf}), or to Elliptic Integrals
if it is half-integer (Section \ref{sec.elli}).

\subsection{Cases of Elementary Functions}\label{sec.elf}

If $\beta=0,1,2,..$ and $k>0$, we substitute $x=\tau^2$
in (\ref{eq.Iofab}),
\begin{equation}
I_{\beta,k}
=\frac{e^{2\beta}}{2(1-c^2e^2)^{k+1/2}}\int_0 dx \frac{(1-x)^{k-1}(a^2-x)^\beta}{\sqrt{x}(b^2-x)^{k+1/2}}
,
\end{equation}
equivalent to the computation of
\begin{equation}
J_{\beta,k}(x)\equiv
\int dx \frac{(1-x)^{k-1}(a^2-x)^\beta}{\sqrt{x}(b^2-x)^{k+1/2}}
.
\label{eq.Jofx}
\end{equation}
Partial integration of this integral
generates the recurrence
\begin{eqnarray}
\left[1+2\beta+\frac{2(k-1)b^2+2k+1}{b^2-1}\right]
J_{\beta,k}
=
2\sqrt{x} \frac{(1-x)^{k-1}(a^2-x)^\beta}{(b^2-x)^{k+1/2}}
+
2\frac{k-1}{b^2-1}J_{\beta,k-1}
\nonumber
\\
+
2a^2\beta J_{\beta-1,k}
+
(2k+1)\frac{b^2}{b^2-1}J_{\beta,k+1}
.
\label{eq.Jrec}
\end{eqnarray}
This allows to build the entire table of $J_{\beta,k}$ from
a list of $J_{\beta,0}$, $J_{\beta,1}$ and $J_{0,k}$.

The substitution
\begin{equation}
z=\frac{1}{b^2-x},\quad  x=b^2-\frac{1}{z}
\end{equation}
and binomial expansion establish
\begin{eqnarray}
J_{\beta,k}
&=& 
\int dz \frac{[(1-b^2)z+1]^{k-1}}{\sqrt{b^2z-1}}
\frac{ [(a^2-b^2)z+1]^\beta}{z^\beta}
\label{eq.Jkbz}
\\
&=&
\sum_{s=0}^{k-1} \binom{k-1}{s} 
\sum_{m=0}^\beta \binom{\beta}{m}
(1-b^2)^s(a^2-b^2)^m
\int dz \frac{z^{s+m-\beta}}{\sqrt{b^2z-1}},\quad k\ge 1,\quad \beta\ge 0.
\label{eq.Jkb}
\end{eqnarray}

\begin{itemize}
\item
The special case $s+m-\beta=0$ is covered by
\begin{equation}
\int \frac{dz}{\sqrt{cz-1}}
=\frac{2}{c}\sqrt{cz-1}.
\end{equation}
\item
The cases $s+m-\beta <0$ are handled by \cite[2.245.2]{GR}
\begin{equation}
\int \frac{dz}{z^{t+1}\sqrt{cz-1}}
=
c^t
\frac{\Gamma(t+1/2)}{\Gamma(t+1)}
\left[
\sqrt{cz-1}
\sum_{l=0}^{t-1}
\frac{\Gamma(l+1)}{\Gamma(3/2+l)}
\frac{1} {(cz)^{l+1}}
+2\frac{1}{\sqrt{\pi}} \arctan\sqrt{cz-1}
\right],\quad t\ge 0.
\end{equation}
The sum evaluates to zero if the upper limit is smaller than the lower limit.
\item
The cases $s+m-\beta\ge 0$ are solved by \cite[2.222]{GR}
\begin{equation}
\int \frac{z^t}{\sqrt{cz-1}} dz = \frac{2\sqrt{cz-1}}{c^{t+1}}\sum_{l=0}^t \binom{t}{l}\frac{ (cz-1)^l}{2l+1}
,\quad t\ge 0
.
\label{eq.intzt}
\end{equation}
\end{itemize}
This has effectively written (\ref{eq.Jkb}) as triple
sums. [In numerical practice,
the integral in (\ref{eq.Jkb}) is placed into a look-up table for
the $\beta+k$ different values of the exponent $s+m-\beta$.]
The case $\beta=0$ appears as a double sum, but is recast into a single sum
by resummation of the $l$-sum in (\ref{eq.intzt}) and the
$s$-sum in (\ref{eq.Jkb}):

\begin{eqnarray}
J_{0,k}&=&
 2
\frac{1}{b^{2k}}
\sqrt{b^2z-1}
\sum_{l=0}^{k-1}
\binom{k-1}{l}(1-b^2)^l
\frac{(b^2z-1)^l}{2l+1}
\\
&=&
2\frac{\sqrt{b^2z-1}}{b^{2k}}
{}_2F_1\left(\frac{1}{2},1-k;\frac{3}{2};(b^2-1)(b^2z-1)\right)
=
2\frac{\sqrt{x}}{b^{2k}\sqrt{b^2-x}}
{}_2F_1\left(\frac{1}{2},1-k;\frac{3}{2};\frac{(1-b^2)x}{b^2-x}\right)
.
\end{eqnarray}

The remaining part of Section \ref{sec.elf} considers
the case $k=0$ which is not available from (\ref{eq.Jkb}):

\begin{eqnarray}
I_{\beta,0}
&=& \frac{e^{2\beta}}{2\sqrt{1-c^2e^2}}
\int d\tau \frac{(a^2-\tau^2)^\beta}{(1-\tau^2)\sqrt{b^2-\tau^2}}
\nonumber \\
&=& \frac{e^{2\beta}}{2\sqrt{1-c^2e^2}}
\sum_{l=0}^\beta \binom{\beta}{l}
(a^2-1)^{\beta-l}
\int d\tau
\frac{(1-\tau^2)^{l-1}}{\sqrt{b^2-\tau^2}}
\nonumber \\
&=&
\frac{e^{2\beta}}{2\sqrt{1-c^2e^2}}
(a^2-1)^\beta
\left[
\frac{1}{2\sqrt{1-b^2}}
\arccos\frac{b^2(1+\tau^2)_2\tau^2}{b^2(1-\tau^2)}
+
\sum_{l=1}^\beta \binom{\beta}{l}
\frac{1}{(a^2-1)^l}
\int d\tau
\frac{(1-\tau^2)^{l-1}}{\sqrt{b^2-\tau^2}}
\right]
.
\label{eq.Ibetazero}
\end{eqnarray}
[The split value at $l=0$ is derived substituting $y=1/(1-\tau^2)$,
then using \cite[2.261]{GR}. The constant of integration is chosen to
realize the limit $I_{0,0}\stackrel{\tau\to 0}{\longrightarrow} 0$.]

This is a superposition of
\begin{equation}
\int d\tau \frac{(1-\tau^2)^{l-1}}{\sqrt{b^2-\tau^2}}
=\sum_{i=0}^{l-1}\binom{l-1}{i} (1-b^2)^{l-1-i}\int d\tau (b^2-\tau^2)^{i-1/2}
=\sum_{i=0}^{l-1}\binom{l-1}{i} (1-b^2)^{l-1-i} B_i,
\label{eq.Bisum}
\end{equation}
where we have defined
\begin{equation}
B_i \equiv \int d\tau (b^2-\tau^2)^{i-1/2};\quad i=0,1,2,.\dots
\end{equation}
Starting from
\begin{eqnarray}
B_0 
&=& \arcsin\frac{\tau}{b}
,
\end{eqnarray}
the recurrence \cite[2.260.2]{GR}
\begin{equation}
B_i=\frac{\tau}{2i}(b^2-\tau^2)^{i-1/2}+\left(1-\frac{1}{2i}\right)b^2B_{i-1}
\end{equation}
allows to calculate (\ref{eq.Bisum}) and eventually (\ref{eq.Ibetazero}).

\subsection{Elliptic case}\label{sec.elli}

This section looks at (\ref{eq.Iofab}) for the cases $\beta=-\frac{1}{2}$, $\frac{1}{2}$, $
\frac{3}{2}$, $\frac{5}{2}\ldots$. Let
\begin{equation}
\bar J_{\beta,k}(\tau)\equiv
2\int d\tau \frac{(1-\tau^2)^{k-1}(a^2-\tau^2)^\beta}{(b^2-\tau^2)^{k+1/2}}
.
\label{eq.Jbar}
\end{equation}
The factor 2 in the definition is chosen to maintain the ``alignment''
$\bar J_{\beta,k}(\tau)=J_{\beta,k}(\tau^2)$.
Partial integration offers the recurrence
\begin{eqnarray}
\left[1+2\beta+\frac{2(k-1)b^2+2k+1}{b^2-1}\right]
\bar J_{\beta,k}
=
2\tau \frac{(1-\tau^2)^{k-1}(a^2-\tau^2)^\beta}{(b^2-\tau^2)^{k+1/2}}
+
2\frac{k-1}{b^2-1}\bar J_{\beta,k-1}
\nonumber
\\
+
2a^2\beta \bar J_{\beta-1,k}
+
(2k+1)\frac{b^2}{b^2-1}\bar J_{\beta,k+1}
,
\end{eqnarray}
which is the same as (\ref{eq.Jrec}).

For $k\ge 1$, binomial expansion of the numerator in (\ref{eq.Jbar})
proposes
\begin{equation}
\bar J_{\beta,k}
=
2
\sum_{s=0}^{k-1}
\sum_{m=0}^{\beta+1/2}
\binom{k-1}{s} (1-b^2)^s
\binom{\beta+1/2}{m} (a^2-b^2)^m
\int 
\frac{d\tau}{\sqrt{(a^2-\tau^2)(b^2-\tau^2)}(b^2-\tau^2)^{1+s-(\beta+1/2-m)}}
.
\end{equation}
These elliptic integrals with integer exponents
$v\equiv 1+s-(\beta+1/2-m)\ge 0$ are \cite[219.06]{Byrd}
\begin{equation}
\int d\tau 
\frac{1}{\sqrt{(a^2-\tau^2)(b^2-\tau^2)}(b^2-\tau^2)^v}
=\frac{1}{ab^{2v}}D_{2v}
,
\end{equation}
where \cite[313.05]{Byrd}\cite[3.158.15]{GR}
\begin{eqnarray}
D_0 &=&
F(\xi,b/a)
,
\\
D_2 &=&
F(\xi,b/a)
+
\frac{a^2}{a^2-b^2}
\left[
\frac{\tau}{a}\frac{\sqrt{a^2-\tau^2}}{\sqrt{b^2-\tau^2}}
-E(\xi,b/a)
\right]
,
\end{eqnarray}
with $\sin\xi\equiv \tau/b$.
[In this equation and where used with two arguments,
(\ref{eq.A2}), (\ref{eq.I32zero}) and the last equation in the appendix, $E$ is the incomplete
Elliptic Integral of the second kind,  elsewhere the shorthand (\ref{eq.TEdef}).]
The recurrence extending these two values both sides
to larger or negative $v$ is \cite[313.05]{Byrd}
\begin{equation}
(2v+1)(a^2-b^2)D_{2v+2} = (2v-1)b^2D_{2v-2}+2v(a^2-2b^2)D_{2v}+a \sqrt{a^2-\tau^2} \frac{\tan\xi}{\cos^{2v}\xi}
.
\end{equation}

Similar to the exception in Section \ref{sec.elf}, the case $k=0$ is not covered
by the expansion above and established individually:
\begin{eqnarray}
\bar J_{\beta,0}
&=&
2 \int d\tau \frac{(a^2-\tau^2)^\beta}{(1-\tau^2)\sqrt{b^2-\tau^2}}
\\
&=&
2 \int d\tau \frac{(a^2-\tau^2)^{\beta+1/2}}{(1-\tau^2)\sqrt{(a^2-\tau^2)(b^2-\tau^2)}}
\\
&=&
2
\sum_{m=0}^{\beta+1/2}
\binom{\beta+1/2}{m}
(a^2-1)^{\beta+1/2-m}
\int d\tau \frac{(1-\tau^2)^{m-1}}{\sqrt{(a^2-\tau^2)(b^2-\tau^2)}}
.
\end{eqnarray}
The term $m=0$
is an Elliptic Integral of the third kind
\cite[219.02]{Byrd}\cite[3.157.7]{GR}:
\begin{equation}
\int d\tau \frac{1}{(1-\tau^2)\sqrt{(a^2-\tau^2)(b^2-\tau^2)}}
=
\frac{1}{a}\Pi(\xi,b^2,b/a)
.
\end{equation}
Because the case $k=s=0$ with $\beta=-1/2$ is the only contribution to 
(\ref{eq.Ialphser}) and (\ref{eq.Ipows}) if the geodetic line is
on the surface of the ellipsoid ($h=0$), this is the only value relevant
to the integral (\ref{eq.I}) for this ``classic'' case.

The terms $m\ge 1$ are delegated
to \cite[219.05]{Byrd} by binomial expansion of the numerator,
\begin{equation}
\int d\tau \frac{(1-\tau^2)^{m-1}}{\sqrt{(a^2-\tau^2)(b^2-\tau^2)}}
=
\sum_{l=0}^{m-1} \binom{m-1}{l}
\frac{(-b^2)^l}{a}A_{2l}
.
\end{equation}
Starting from \cite[310.02,310.05]{Byrd}
\begin{eqnarray}
A_0&=& F(\xi,b/a);
\\
b^2A_2&=& a^2[F(\xi,b/a)-E(\xi,b/a)],
\label{eq.A2}
\end{eqnarray}
more values of $b^{2l}A_{2l}$
may be generated from the recurrence \cite[310.05]{Byrd}
\begin{equation}
(2l+1)b^{2l+2}A_{2l+2}
= 
a\sqrt{b^2-\tau^2}\sqrt{a^2-\tau^2} +2l(a^2+b^2)b^{2l}A_{2l}+(1-2l)a^2b^{2l}A_{2l-2}
;\quad l\ge -1
.
\end{equation}

\section{Line Distance Integral}\label{sec.sint}
\subsection{Reduction to the Angular Coupling Integral}

The formula for the distance $s$ along the geodetic line
is given by \cite{MatharArxiv0711}
\begin{eqnarray}
s&=&
\int d\tau
\frac{
h(1-e^2\tau^2)+N(1-e^2)
}{
(1-e^2\tau^2)
\sqrt{ 1-\tau^2-\frac{c^2}{(N+h)^2} }
}
\label{eq.sIntgr}
\\
&=&
\rho_e \int d\tau
\frac{
(1+\frac{h}{\rho_e}E^{1/2})\frac{h}{\rho_e}E^{3/2}
+1-e^2+\frac{h}{\rho_e}(1-e^2)E^{1/2}
}{
E^{3/2}
\sqrt{(1+\frac{h}{\rho_e}E^{1/2})^2 T-\frac{c^2}{\rho_e^2}E} }
\\
&=&
\rho_e
\left[
\frac{h}{\rho_e}
S_0(\frac{h}{\rho_e},\frac{c}{\rho_e})
+
(\frac{h}{\rho_e})^2
S_{1/2}(\frac{h}{\rho_e},\frac{c}{\rho_e})
+
(1-e^2)
S_{-3/2}(\frac{h}{\rho_e},\frac{c}{\rho_e})
+
\frac{h}{\rho_e}(1-e^2)
S_{-1}(\frac{h}{\rho_e},\frac{c}{\rho_e})
\right]
,
\label{eq.Alist}
\end{eqnarray}
calling a group of integrals
\begin{equation}
S_{\alpha}(h,c)\equiv
\int d\tau
\frac{E^\alpha}{ \sqrt{T +2hE^{1/2} T +h^2 ET -c^2E} }
\end{equation}
with two dimensionless scaled parameters and one characteristic exponent $\alpha$.
Binomial expansion of the square root provides a power series in $h$,
\begin{equation}
S_{\alpha}(h,c)
=
\sum_{s=0}^\infty
h^s
\sum_{k=\lceil s/2\rceil}^s
\binom{-1/2}{k}
\binom{k}{s-k}
2^{2k-s}\int d\tau \frac{E^{\alpha+s/2}T^k}{(T-c^2E)^{k+1/2}}
.
\label{eq.Sofs}
\end{equation}
The integral is similar to
(\ref{eq.Ibetakdef});
the difference is a factor $T$
plus the request from (\ref{eq.Alist}) to evaluate the cases $\alpha=-3/2$ and $\alpha=-1$.
The factor $T$ is distributed with the aid of
\begin{equation}
T=1-\frac{1}{e^2}+\frac{1}{e^2}E
.
\end{equation}
Recalling (\ref{eq.adef}), this maps (\ref{eq.Sofs}) on the integrals (\ref{eq.Ibetakdef})
\begin{equation}
\int d\tau \frac{E^{\alpha+s/2}T^k}{(T-c^2E)^{k+1/2}}
=
(1-a)I_{\alpha+s/2,k}
+
a I_{\alpha+1+s/2,k}.
\label{eq.SofI}
\end{equation}

\subsection{Special Values}
To carry out the right hand side of (\ref{eq.SofI}), the
only aspect not yet covered by Section \ref{sec.llc} is
to implement $I_{-3/2,k}$ and $I_{-1,k}$.
Furthermore, $k$ is only required for the restricted
range of $s$ seen in the summation (\ref{eq.Sofs}), which reduces the ``new'' cases further
to
\begin{itemize}
\item
$I_{-3/2,0}$ from $(s,k,\alpha)=(0,0,-3/2)$,
\item
$I_{-1,0}$ from $(s,k,\alpha)=(0,0,-1)$,
\item
and $I_{-1,1}$ from $(s,k,\alpha)=(1,1,-3/2)$,
\end{itemize}
because otherwise the first index of $I_{\beta,k}$ is $\beta\ge -1/2$,
already treated in Section \ref{sec.llc}.
Turning to
\begin{equation}
I_{-3/2,0}
=
\frac{e^{-3}}{(1-c^2e^2)^{1/2}}
\int d\tau \frac{1}{(1-\tau^2)(a^2-\tau^2)^{3/2}(b^2-\tau^2)^{1/2}}
\end{equation}
proposes partial fraction decomposition
\begin{eqnarray}
\int d\tau \frac{1}{(1-\tau^2)(a^2-\tau^2)^{3/2}(b^2-\tau^2)^{1/2}}
=&&
\frac{1}{a^2-1}\int d\tau \frac{1}{(1-\tau^2)\sqrt{(a^2-\tau^2)(b^2-\tau^2)}}
\nonumber \\
&&
-\frac{1}{a^2-1}\int d\tau \frac{1}{(a^2-\tau^2)\sqrt{(a^2-\tau^2)(b^2-\tau^2}}
,
\label{eq.I32pf}
\end{eqnarray}
and the two integrals on the right hand side are
known \cite[219.02,219.07,315.02]{Byrd}:
\begin{equation}
I_{-3/2,0}=\frac{1}{e(a^2-1)\sqrt{1-c^2e^2}}
\left\{
\Pi(\xi,b^2,b/a)
-
\frac{1}{a^2-b^2}\left[ E(\xi,b/a)-\frac{\tau}{a} \frac{\sqrt{b^2-\tau^2}}{\sqrt{a^2-\tau^2}}\right]
\right\}
.
\label{eq.I32zero}
\end{equation}
Demonstrated in (\ref{eq.Alist}), this is the only $S_\alpha$ value required 
on the surface of the ellipsoid---where $h=0$.

The second remaining case is
\begin{equation}
I_{-1,0}
=
\frac{a}{\sqrt{1-c^2e^2}}
\int d\tau\frac{1}{(1-\tau^2)(a^2-\tau^2)\sqrt{b^2-\tau^2}}
,
\end{equation}
which splits into two partial fractions---equivalent to (\ref{eq.I32pf})---with
known integrals \cite[2.284]{GR}:
\begin{equation}
I_{-1,0}=\frac{a}{\sqrt{1-c^2e^2}(a^2-1)}
\left[
\frac{1}{\sqrt{1-b^2}}\arctan\left(\tau \sqrt{\frac{1-b^2}{b^2-\tau^2}}\right)
-
\frac{1}{a\sqrt{a^2-b^2}}\arctan\left(\frac{\tau}{a} \sqrt{\frac{a^2-b^2}{b^2-\tau^2}}\right)
\right]
.
\end{equation}

The third remaining case is
\begin{equation}
I_{-1,1}
=
\frac{a}{(1-c^2e^2)^{3/2}}
\int \frac{d\tau}{(a^2-\tau^2)(b^2-\tau^2)^{3/2}}
,
\end{equation}
with partial fractions
\begin{equation}
\int \frac{d\tau}{(a^2-\tau^2)(b^2-\tau^2)^{3/2}}
=
\frac{1}{b^2-a^2}\int \frac{d\tau}{(a^2-\tau^2)\sqrt{b^2-\tau^2}}
-
\frac{1}{b^2-a^2}\int \frac{d\tau}{(b^2-\tau^2)^{3/2}}
.
\end{equation}
The first integral on the right hand side is the same as met
while calculating $I_{-1,0}$.
The second is also known \cite[2.271.5]{GR}, to yield
\begin{equation}
I_{-1,1}
=
\frac{1}{(1-c^2e^2)^{3/2}(a^2-b^2)}
\left[
\frac{a}{b^2}\frac{\tau}{\sqrt{b^2-\tau^2}}
-
\frac{1}{\sqrt{a^2-b^2}}
\arctan\left(\frac{\tau}{a}\sqrt{\frac{a^2-b^2}{b^2-\tau^2}}\right)
\right]
.
\end{equation}

\section{Inverse Function}\label{sec.inv}

In sections \ref{sec.llc} and \ref{sec.sint}
the value of the parameter $c$ was regarded as known.
The ``inverse'' problem of geodesy, on the other hand,
is finding $c$ assuming the value of the integral, i.e., $\Delta \lambda$,
and of its limits are given. If Newton methods are employed to this
problem, they also call for computation of 
the power series of $\partial I/\partial c$, which is addressed as follows:

The derivative of (\ref{eq.IofIalph}) is
\begin{equation}
\partial_c I=\frac{1}{\rho_e}I+\frac{c}{\rho_e}\left[
\frac{h}{\rho_e^2}\partial_c I_1(h,c)+\frac{1}{\rho_e}(1-e^2)\partial_c I_{-1/2}(h,c)
\right]
\end{equation}
by the product rule.
The derivative of (\ref{eq.Ialphser}) is given by the chain rule:
\begin{equation}
\partial_c I_\alpha(h,c)
= c\sum_{s=0}^\infty (-h)^s \sum_{k=0}^s (2k+1)\kappa_{s,k}
\int d\tau \frac{ E^{1+\alpha+s/2}}{(T-c^2E)^{k+3/2}} T^{k-1}
.
\end{equation}
The exponent of $T$ in this integrand is deficient by 1 compared
to the definition of $I_{\beta,k}$ in (\ref{eq.Ibetakdef})---whereas
it is abundant by 1 in (\ref{eq.Sofs}). Still, new functionality is
not required, because
\begin{equation}
\int d\tau \frac{ E^{1+\alpha+s/2}}{(T-c^2E)^{k+3/2}} T^{k-1}
=
\int d\tau \frac{1}{T(T-c^2E)}\frac{ E^{1+\alpha+s/2}}{(T-c^2E)^{k+1/2}} T^k
\end{equation}
\begin{equation}
=
\int d\tau \frac{1}{c^2E}\left[\frac{1}{T-c^2E}-\frac{1}{T}\right]
\frac{ E^{1+\alpha+s/2}}{(T-c^2E)^{k+1/2}} T^k
=
\frac{1}{c^2} I_{\alpha+s/2,k+1}
-
\frac{1}{c^2} I_{\alpha+s/2,k}
\end{equation}
converts to integrals
already discussed in Section \ref{sec.llc}.

\section{Summary}
Given the altitude, a directional parameter and a starting position,
the task of finding the trajectory of the geodetic line in 3-dimensional geodetic coordinates
turns into the evaluation of an integral which emerges from the solution of a 
differential equation which couples
longitude and latitude.
The coefficients of a systematic expansion of this integral in a power series
of the altitude (scaled by the equatorial radius) have been reduced to multiple
sums over elementary functions and incomplete Elliptic Integrals;
the geodetic line on the surface of the ellipsoid is a specific and the simplest case.
Auxiliary integrals developed along this path recur
if the same expansion strategy is applied to other integrals
related to the geodetic line.

\appendix

\section{Table Errata}
Errata to the 1981 edition of the Tables of Sums, Products and Integrals \cite{GR} relevant to this work are:
\begin{itemize}
\item
2.284: Preserve the sign of $c$ on the right hand side by moving it out of the square root:
\[
\int \frac{Ax+B}{(p+R)\sqrt{R}}dx
=
\frac{A}{c}I_1+\frac{2Bc-Ab}{c\sqrt{p[b^2-4(a+p)c]}}I_2.
\]
\item 2.245.2:
Alternate signs of both $b$ on the right hand side:
\begin{eqnarray*}
\int\frac{z^m dx}{t^n\sqrt{z}}&=& -z^m\sqrt{z}\big\{\frac{1}{(n-1)\Delta}\frac{1}{t^{n-1}}
\\
&&+\sum_{k=2}^{n-1}\frac{(2n-2m-3)(2n-2m-5)\cdots(2n-2m-2k+1)(-b)^{k-1}}
{2^{k-1}(n-1)(n-2)\cdots(n-k)\Delta^k}\,\frac{1}{t^{n-k}}\big\}
\\
&&
-\frac{(2n-2m-3)(2n-2m-5)\cdots(-2m+3)(-2m+1)(-b)^{n-1}}{2^{n-1}(n-1)!\Delta^n}\int \frac{z^m dx}{t\sqrt{z}}.
\end{eqnarray*}
\item 3.158.15:
Remove a $b$ in a numerator of the right hand side:
\begin{eqnarray*}
\int_0^u \frac{dx}{\sqrt{(a^2-x^2)(b^2-x^2)^3}}
&=& \frac{1}{ab^2}F(\eta,t)-\frac{1}{b^2(a^2-b^2)}
\\
&&
\times
\left\{aE(\eta,t)-u\sqrt{\frac{a^2-u^2}{b^2-u^2}}\right\},\quad [a>b>u>0].
\end{eqnarray*}
\end{itemize}

\bibliographystyle{apsrmp}

\bibliography{all}

\end{document}